\documentclass [11pt]{amsart} 
\usepackage{epic,eepic,latexsym, amssymb, amscd, amsfonts}

\input xy
\xyoption {all}

\def \quot {{{\text{Quot}}_{d} ({\mathcal O}^{r+k}, k, C)}}
\def \gbar {{\bar g}}

\def \jac {{\text{Jac}}}
\def \mor {{{\text{Mor}}_{d}(C, G(k, r+k))}}
\def \mors {{{\text{Mor}}_{d}^{\text{stable}}(C, G(k, r+k))}}

\newtheorem{theorem}{Theorem} 
\newtheorem {lemma}{Lemma}
 
\newtheorem {corollary}{Corollary} 
\newtheorem {proposition}{Proposition}

\theoremstyle{definition}
\newtheorem{remark}{Remark} 

\theoremstyle {definition}

\begin{document}

\title[The rank-level duality for non-abelian theta functions]{The rank-level duality for non-abelian theta functions}
\author {Alina Marian}
\address {Department of Mathematics}
\address {Yale University}
\email {alina.marian@yale.edu}
\author {Dragos Oprea}
\address {Department of Mathematics}
\address {Stanford University}
\email {oprea@math.stanford.edu}
\date{} 

\begin {abstract} We prove that the strange duality conjecture of Beauville-Donagi-Tu holds for all 
curves. We establish first a more extended rank-level duality, interesting in its own right, from which the standard 
rank-level duality follows by restriction.
\end{abstract} \maketitle

\section{Introduction}

Let $C$ be a smooth complex projective curve of genus $g\geq 1$. The Jacobian variety $\jac^{0}(C)$, parametrizing degree $0$ line
bundles on the curve $C$, comes equipped with a Theta divisor, unique up to translation. For any line bundle $M$ on $C$ of degree
$\gbar = g-1$, the divisor $\Theta_{1,M}$ is given set theoretically by $$\Theta_{1,M}=\left\{E \text { degree } 0 \text{ line bundle, such that } h^{0}(E\otimes M)=h^{1}(E\otimes M)\neq 0\right\}.$$ The sections of higher tensor powers $\mathcal 
O(l\Theta_{1,M})$ of the associated Theta 
line bundle are the well-studied theta functions of level $l$.

The non-abelian version of this setup replaces the Jacobian by the moduli space $U(k,d)$ of semistable bundles of rank $k$ 
and degree $d$. In the particular case when $d=k\gbar$, the moduli space $U(k, k\gbar)$ carries a canonical Theta divisor 
$\Theta_k$. On the stable part of the moduli space, $\Theta_k$ is supported on the jumping locus $$\Theta_k= \left\{E\in U(k,k\gbar) \text { stable, such that } h^{0}(E)=h^{1}(E)\neq 0\right\}.$$ Since $U(k, k\gbar)$ is locally factorial \cite {DN}, this divisor corresponds to a line bundle also 
denoted 
$\Theta_{k}$. The context will unambiguously indicate throughout the paper whether we refer to the divisor or to the 
associated line bundle. To construct Theta divisors for arbitrary degree, consider first a vector bundle $M$ of rank $r$ on 
$C$ with slope $$\mu=\gbar-\frac{d}{k}.$$ This numerical condition is equivalent to $$\chi(E\otimes M)=0 \text { for all } 
E\in U(k,d).$$ Letting $$\tau_{M}: U (k,d) \longrightarrow U(kr, kr\gbar), \, \, \, \, E \mapsto E \otimes M$$ be the 
tensor product map, we define the line bundle $$\Theta_{k, M} = \tau_M^{\star} \Theta_{kr}.$$ $\Theta_{k,M}$ has a 
canonical pullback section which, on the stable part of the moduli space, vanishes on the jumping locus \begin{equation} 
\label{thetam} \Theta_{k,M}=\left\{E\in U(k,d)\text { stable, such that } h^0(E\otimes M)=h^{1}(E\otimes M)\neq 0\right\}. 
\end{equation} For a {\it generic} choice of $M$, the vanishing locus \eqref{thetam} is a divisor on $U(k,d)$ \cite {RT}.

The same construction produces Theta divisors on the moduli space $SU(k, \Lambda)$ of semistable bundles with fixed 
determinant $\Lambda$. The associated line bundles only depend on the rank of the bundle $M$ used in the definition. When 
$M$ has minimal possible rank, the corresponding line bundle $\theta_{k}$ is an ample generator for the Picard group of 
$SU(k, \Lambda)$ \cite {DN}.

The generalized ({\it i.e.}, non-abelian) theta functions (of rank $k$ and level $l$) are the global sections of $l$-tensor
powers $\Theta_k^l$ of the Theta bundles $\Theta_{k}$. They are conjectured to be related by a ``strange duality''
isomorphism which interchanges the rank and the level. The duality is defined via the tensor product map
$$\tau:SU(r,\mathcal O)\times U(k, k\gbar) \to U(kr, kr\gbar),$$ given by $$(E,F)\to E\otimes F.$$ It is easy to prove
using the see-saw principle that $$\tau^{\star} \Theta_{kr}\cong\theta_r^k\boxtimes \Theta_{k}^r.$$ The canonical section
of the line bundle $\Theta_{kr}$ determines by pullback an element of
$$H^{0}(SU(r,\mathcal O)\times U(k, k\gbar) , \tau^{\star}\Theta_{kr})=H^{0}(SU(r,\mathcal O), \theta_r^{k})\otimes
H^{0}(U(k, k\gbar), \Theta_{k}^{r}).$$ In turn, this induces the strange duality morphism $\mathsf {SD}$, well defined up
to scalars, \begin{equation}\label{sd}\mathsf{SD}: H^{0}(SU(r,\mathcal O), \theta_r^{k})^{\vee} \to H^{0}(U(k, k\gbar),
\Theta_{k}^{r}).\end {equation}

It has been conjectured that $\mathsf {SD}$ is an isomorphism.  In the mathematics literature, the conjecture is stated in
\cite{Bea}, \cite{DT}. Its first supporting fact is that the dimensions of the two spaces, given by the Verlinde formula,
can be explicitly seen to be the same. Proofs of the duality in particular cases can be found in \cite {Bea1}, \cite{BNR},
\cite{l}, \cite{vGP}, \cite {P}. Most notably, P. Belkale \cite{Bel} recently established the isomorphism \eqref{sd} for
generic curves. On the other hand, the physics literature - for instance \cite {NS} and the references therein - places the
discussion in the context of conformal blocks and the Wess-Zumino-Witten model. Indeed, the representation theory of affine
Lie algebras seems to be a natural home for such rank-level dualities, and statements of related flavor have been
formulated in \cite {F}.

In this paper, we prove first a more extended rank-level duality, which holds for all curves, and relates nonabelian theta
functions on moduli spaces of bundles with {\it varying} determinant. We will show

\begin{theorem}\label{main2} Let $L$ and $M$ be any two line bundles of degree $\gbar$ on $C$, and let $K$ be the canonical bundle of $C$. Then there exists a rank-level duality isomorphism between the spaces of global sections 
\begin{equation}
\label{sd1alt}
\mathsf {D}:H^{0}\left (U(r,0),
\Theta_{r,M}^{k}\otimes {\det}^{\star} \Theta_{1,L}\right)^{\vee} \to H^{0} \left (U(k,0), \Theta_{k,M}^{r}\otimes 
{\det}^{\star}
\Theta_{1,L^{-1}\otimes K} \right ).
\end{equation}
Here $$\det: U(r,0)\to \jac^0(C) \text { and } \det: U(k, 0)\to 
\jac^{0}(C)$$
are the morphisms
taking bundles to their determinants.  
\end{theorem}

The theorem, in conjunction with a restriction argument involving the geometry of the classical theta functions on the 
Jacobian $\jac^{0} (C)$, gives the following important strengthening of the main result of \cite{Bel}.

\begin{theorem} \label{main} The strange duality morphism $\mathsf {SD}$ is an isomorphism for all curves $C$. \end
{theorem}

The morphism $\mathsf {D}$ will be defined in Section 2.  We briefly describe it here in the simplest case $ r = k = 1.$ Let
$$\rho:\jac^{0}(C)\times \jac^{0}(C)\to \jac^{0}(C)\times \jac^{0}(C)$$ be the isogeny given by $$(A, B)\to (A\otimes B, A^{-1}\otimes B),$$ and let $$(-1):\jac^{0}(C)\to \jac^{0}(C)$$ denote the multiplication by $(-1)$ in the Jacobian. An
immediate application of the see-saw principle proves the isomorphism $$\rho^{\star}\left(\Theta_{1, L}\boxtimes
\Theta_{1,M}\right)\cong\left(\Theta_{1,L}\otimes \Theta_{1,M}\right) \boxtimes \left((-1)^{\star}\Theta_{1,L}\otimes \Theta_{1,M}\right).$$ Then, the canonical pullback section of $\rho^{\star}\left(\Theta_{1,
L}\boxtimes \Theta_{1,M}\right)$ determines the duality morphism $${\mathsf {D}}: {H^0 \left (\jac^0(C), \Theta_{1,L}\otimes \Theta_{1,M} \right
)}^{\vee} \longrightarrow H^0 \left (\jac^0(C), (-1)^{\star}\Theta_{1,L}\otimes \Theta_{1,M}\right ).$$ Let us remark that $\mathsf {D}$ can be explicitly
diagonalized by a particular choice of a Theta-basis. While the general theory applies to arbitrary abelian varieties, hence to any Jacobian, we
will assume here, for simplicity of notation, that $C=E_{\tau}$ is an elliptic curve with complex modulus $\tau$, and that $L$ and $M$ coincide
and correspond to the point $\frac{\tau+1}{2}$ of $E_{\tau}$. In this case, $\Theta = \Theta_{1, L}$ has the canonical section
$$\theta_{0,0}^{\tau}(z)=\sum_{n\in \mathbb Z} \exp(\pi i n^{2}\tau+2\pi i n z).$$ The above discussion shows that $\mathsf {D}$ is determined
by the section $$\rho^{\star}\left(\theta^{\tau}_{0,0}(z)\boxtimes
\theta_{0,0}^{\tau}(w)\right)=\theta^{\tau}_{0,0}(z+w)\theta^{\tau}_{0,0}(w-z)\in H^{0}(\Theta^2)\otimes H^{0}(\Theta^2).$$ Different {\it
explicit} tensor product decompositions of this section in terms of various bases of level $2$ theta functions are known as addition formulas
(\cite {LB}, page 208). For instance, if we let $$\theta^{\tau}_{\frac{1}{2},\frac{1}{2}}(z)=\exp\left(\frac{\pi i \tau}{4}+\pi i
\left(z+\frac{1}{2}\right)\right)\theta_{0,0}^{\tau}\left(z+\frac{\tau+1}{2}\right)$$ be the scaled translate of $\theta^{\tau}_{0,0}(z)$, then
$$\left(\theta^{\tau}_{0,0}\right)^2, \left(\theta^{\tau}_{\frac{1}{2},\frac{1}{2}}\right)^2\in H^{0}(\Theta^{2})$$ form a basis for the level
$2$ theta functions. The corresponding addition formula is \cite {M}:  
$$\theta^{\tau}_{0,0}(z+w)\theta^{\tau}_{0,0}(w-z)\theta^{\tau}_{0,0}(0)^2=\theta^{\tau}_{0,0}(z)^2\theta^{\tau}_{0,0}(w)^2+\theta^{\tau}_{\frac{1}{2},
\frac{1}{2}}(z)^2\theta^{\tau}_{\frac{1}{2}, \frac{1}{2}}(w)^2.$$ Therefore, $\mathsf {D}$ is a multiple of the identity in the above theta
basis. The reader is probably aware of the many companion addition formulas \cite {M}; they determine different matrix representations for the
duality map. \vskip.1in

The setup of Theorems 1 and 2 easily generalizes to the case of arbitrary degree. Specifically, we let $d$ and $r$ be coprime 
integers, and
$h$, $k$ be any two non-negative integers. In addition, we fix a reference bundle $Q$ of degree $d$ and rank $r$. Then,
imitating the construction leading to \eqref{sd}, we obtain a morphism
\begin{equation}\label{sdd}\mathsf{SD}:H^{0}\left(SU\left(hr, (\det Q)^{h}\right), \theta^k_{hr}\right)^{\vee}\to
H^{0}\left(U\left(kr, k(r\bar{g} - d)\right), \Theta_{kr,Q}^{h}\right).\end{equation} The explicit definition of $\mathsf {SD}$ is 
contained in the
last section of this paper. There, we describe first an extended duality morphism, showing

\begin{theorem} \label{arbdegree} There is a rank-level duality isomorphism $\mathsf {D}$ between the spaces of global sections
$$H^{0}\left(U(hr, hd), \Theta_{hr, M\otimes Q^{\vee}}^{k}\otimes {\det}^{\star} \Theta_{1, L\otimes {\left(\det
Q\right)}^{-h}}\right)^{\vee}, \text { and }$$ $$H^{0}\left(U(kr, -kd), \Theta_{kr, M\otimes Q}^{h}\otimes
{\det}^{\star}\Theta_{1, L^{-1}\otimes K \otimes \left(\det Q\right)^{k}}\right).$$ \end{theorem}

As a corollary, by restriction to the moduli space with fixed determinant, we confirm the following conjecture of Donagi-Tu \cite {DT}.

\begin{theorem} \label{arbsd}The rank-level duality map \eqref{sdd} is an isomorphism.

\end{theorem}\vskip.1in Let us briefly outline the structure of our proofs. When the degree is $0$, upon checking equality 
of dimensions for the two spaces of sections involved, the central part of the argument will show that the duality morphism 
$\mathsf {D}$ of Theorem \ref{main2} is surjective. The general strategy is similar to the approach taken in \cite {Bel}. 
There, the analogous question for the {\it original} strange duality morphism is rephrased as a counting problem for 
bundles with certain numerical and cohomological properties. The enumerative problem is first set on a rational nodal 
curve, and its solution is then moved to neighboring smooth curves of genus $g$, hence the requirement that $C$ be generic. 
The main observation of the current paper is that the study of the {\it new} duality morphism $\mathsf {D}$ leads to an 
enumerative problem which can be effectively solved without specializing the curve. For any smooth $C$ of genus $g$ one 
relies instead 
on a strikingly 
simple interpretation of the Verlinde formula providing the dimension of the space of generalized theta functions 
$$H^{0}(U(r,0), \Theta_{r,M}^{k}\otimes {\det}^{\star} \Theta_{1,L}),$$ in terms of counts of maps from $C$ to 
the Grassmannian $G(k, r+k)$. This recasting of the Verlinde formula is due to Witten \cite {witten}, who used physical 
arguments to relate the sigma model of the Grassmannian to Wess-Zumino-Witten theory.  More precisely, the above Verlinde 
number is a top intersection product on the Grothendieck Quot scheme $\quot$ of rank $k$ subsheaves of the trivial bundle 
of rank $r+k$. This intersection product is realized geometrically by points in the Quot scheme, corresponding to stable 
quotients with the numerical and cohomological properties needed to prove the surjectivity of the morphism $\mathsf {D}$. 
That the strange duality morphism \eqref{sd} is an isomorphism will be shown by restricting $\mathsf {D}$ to the moduli 
space of bundles with fixed determinant.

We follow a similar strategy to establish Theorems \ref{arbdegree} and \ref{arbsd}. We rely foremost on an explicit 
elementary expression for the Verlinde formula valid in arbitrary degree. Since this seems to be less known than its degree 
$0$ analogue, found for instance in \cite {Bea}, we decided to record it in the last section of the paper. \vskip.1in

The paper is organized as follows. In the next section, we define the modified duality morphism $\mathsf{D},$ and prove 
that it is an isomorphism, thereby establishing Theorem 1. The argument uses the main enumerative result of Section 3. 
There, we realize the Verlinde number as a count of points in the Quot scheme. In the fourth section we prove Theorem 
2, the original strange duality conjecture, via a general argument on the Jacobian. Finally, the last section discusses the 
strange duality isomorphism in arbitrary degree.

\vskip.1in {\bf Acknowledgements.} 
We are very grateful to Prakash Belkale for his interest in this work, and for
enthusiastic and extensive email correspondence. This paper draws inspiration from his article \cite{Bel}. We thank Mihnea 
Popa for helpful and stimulating conversations on the subject of strange dualities. 

\section{A new rank-level duality map}

In this section we give a precise definition of our new rank-level duality morphism $\mathsf {D}$ announced in the introduction.
We further show that the main enumerative result of Section 3 implies that $\mathsf {D}$ is an isomorphism.\vskip.1in

To begin, fix two line bundles $L$ and $M$ of degree $\gbar$ on $C$. We define the following divisor, denoted $\Delta_{L, M}$, on the product
space $U (r, 0)\times U(k, 0)$. On the stable part of the moduli space, $\Delta_{L,M}$ is given set theoretically as $$\Delta_{L, M} = \{ (E,
F)\text{ such that } h^0 (E \otimes F\otimes M) \neq 0 \text{ or } h^0 ({\det E}^{\vee} \otimes {\det F} \otimes L^{-1}\otimes K) \neq 0 \}.$$
Scheme-theoretically, if we let $$\pi: U(r, 0)\times U(k, 0)\to U(kr, 0)\times \jac^{0}(C)$$ be the morphism $$(E, F)\to (E\otimes F, \det 
E^{\vee}\otimes
\det F),$$ then the line bundle $$\mathcal O(\Delta_{L,M})=\pi^{\star} \left(\Theta_{kr, M}\boxtimes \Theta_{1, L^{-1}\otimes K}\right)$$ has a
natural section which vanishes along the jumping locus $\Delta_{L,M}$ considered above.

\begin{lemma} \label{cross} The following equality holds on $U(r,0)\times U(k,0)$: $$\mathcal O(\Delta_{L, M}) = 
\left(\Theta_{r, M}^k \otimes {\det}^{\star}   
\Theta_{1, L} \right)\boxtimes \left(\Theta_{k,M}^r \otimes {\det}^{\star} \Theta_{1, L^{-1}\otimes K}\right).$$ \end {lemma}

{\it Proof.} The lemma is a consequence of the see-saw principle. We will need the following general fact established in 
\cite{DN}, and which 
will
be invoked several times in this paper.\vskip.08in {\it {\bf Fact. } Assuming $T_1$ and $T_2$ are vector bundles with
$$\mu(T_1)=\mu(T_2)=\gbar, \text{ and }\text {rank } T_1=m\text{ rank } T_2,$$ we have the following isomorphism of Theta bundles on $U(r,0)$:
$$\Theta_{k,T_1}\cong\Theta_{k, T_2}^{m}\otimes {\det}^{\star}(\det T_1\otimes \det T_2^{-m}).$$ Here, we regard the degree $0$ line bundle
$\det T_1\otimes \det T_2^{-m}$ on $C$ as a line bundle on $\jac^{0}(C)$ in the usual way.}\vskip.08in

Fixing $E\in U(r,0)$, we evaluate the restriction $\mathcal O(\Delta_{L, M})$ to
$\{E\}\times U (k, 0)$: \begin{eqnarray}\mathcal O(\Delta_{L,M})\big{|}_{\{E\}\times U(k,0)}&=&\Theta_{k, E\otimes M} \otimes {\det}^{\star} \left(\Theta_{1, \det E^{\vee}\otimes
L^{-1}\otimes K}\right)\cong\nonumber\\&\cong&\left(\Theta_{k,M}^{r}\otimes {\det}^{\star} (\det E)\right) \otimes {\det}^{\star}\left(\Theta_{1, L^{-1}\otimes
K}\otimes \det E^{\vee}\right)\cong\nonumber\\ &\cong&\Theta_{k,M}^{r}\otimes {\det}^{\star}\Theta_{1, L^{-1}\otimes K}.\nonumber\end{eqnarray} 

Next, write $$(-1):\jac^{0}(C)\to \jac^{0}(C)$$ for the multiplication by $(-1)$. As a consequence of (relative) Serre duality, we
have 
$$(-1)^{\star} \Theta_{1, P^{-1}\otimes K}\cong\Theta_{1, P},$$
for all line bundles $P$ of degree $\gbar$. The restriction of
$\mathcal O(\Delta_{L, M})$ to ${U(r,0) \times \{F\} }$ is then seen to be: \begin{eqnarray}\mathcal O(\Delta_{L,M})\big{|}_{U(r,
0)\times \{F\}}&=&\Theta_{r, F\otimes M}\otimes \left({(-1)\circ\det}\right)^{\star}\Theta_{1, \det F \otimes L^{-1}\otimes
K}\cong\nonumber\\&\cong&\left(\Theta_{r,M}^{k}\otimes {\det}^{\star}(\det F)\right)\otimes {\det}^{\star}(-1)^{\star}\Theta_{1,
\det F\otimes L^{-1}\otimes K}\cong\nonumber\\ &\cong&\left(\Theta_{r,M}^{k}\otimes {\det}^{\star}(\det F)\right)\otimes
{\det}^{\star}\Theta_{1, \det F^{\vee} \otimes L}\cong\nonumber\\&\cong&\left(\Theta_{r,M}^{k}\otimes {\det}^{\star}(\det
F)\right)\otimes{\det}^{\star}\left(\Theta_{1, L}\otimes \det F^{\vee}\right)\cong\nonumber\\&\cong&\Theta_{r,M}^{k}\otimes
{\det}^{\star} \Theta_{1, L}.\nonumber\end{eqnarray} This concludes the proof of the lemma. \qed\vskip.1in

Now the natural section $\Delta_{L, M}$ of $\mathcal O(\Delta_{L, M})$ gives the projectively well-defined duality map 
\begin{equation}\label{sd1}\mathsf {D}: H^0 (U(r,
0), \Theta_{r, M}^k \otimes {\det}^{\star} \Theta_{1, L} )^{\vee} \longrightarrow H^0 (U (k, 0), \Theta_{k,M}^r \otimes
{\det}^{\star} \Theta_{1, L^{-1}\otimes K} ),\end{equation} which is the subject of Theorem \ref{main2}. We now give the proof of the theorem, invoking the enumerative conclusions of the next section. \vskip.1in

{\it Proof of Theorem \ref{main2}.} We will note in Section $3$, by means of Corollary \ref{sd2}, that both the domain and the target of the duality map $\mathsf {D}$ have the same dimension. Thus in order to establish Theorem \ref{main2}, it suffices to prove that $\mathsf {D}$ is surjective. 
This will follow easily upon giving a suitable 
enumerative interpretation to the Verlinde number
$$ q = h^0 (U(r, 0), \Theta_{r, M}^k \otimes {\det}^{\star} \Theta_{1, L} ).$$
Indeed, Proposition \ref{pairs} of the next section constructs $q$ pairs  
$(A_i, B_i), \, \, 1 \leq i \leq q$ such that \begin{itemize} \item $A_1, \ldots, A_q$ are rank $r$ degree $0$ stable 
bundles \item $B_1, \ldots, B_q$ are rank $k$ degree $0$ stable  
bundles \end{itemize} and \begin{itemize} \item [(a)] $h^0(A_i \otimes B_j \otimes M) =0$ if and only if $i=j$. \item [(b)] $h^{0}(\det A_i^{\vee}\otimes \det B_i\otimes L^{-1}\otimes K) = 0.$ \end{itemize} The existence of these pairs of bundles implies the surjectivity of the duality map $\mathsf{D}$. 

Indeed, for each $E \in U(r, 0)$, define $$\Delta_E = \{F \in U (k, 0)\, \text{s.t.} \, h^0 (E \otimes F\otimes M) \neq 0 \text { or } h^0 (\det E^{\vee} \otimes \det F \otimes L^{-1}\otimes K) \neq 0\}.$$ It is clear that the
$\Delta_E$s are in the image of $\mathsf {D}$. The two properties (a) and (b) imply that for the pairs $(A_i, B_i), \, \, 1 \leq i \leq q$, we have
$$\Delta_{A_i}(B_j)\neq 0 \text { iff } i=j.$$ It follows that $\Delta_{A_i}, 1\leq i\leq q$ are linearly independent divisors in the image of $\mathsf
{D}$, and for dimension reasons they must span $$H^{0}(U(k,0), \Theta_{k,M}^{r}\otimes {\det}^{\star}\Theta_{1,
L^{-1}\otimes K}).$$\qed

\section{Intersections on the Quot scheme and a formula of Witten}

The main result of this section, Proposition \ref {pairs}, constructs the pairs of vector bundles $(A_1, B_1), \ldots, (A_q, B_q)$
which, as was just noted, render Theorem \ref{main2} evident. These pairs are obtained as points in a suitable Quot scheme; they
geometrically realize an intersection product on Quot whose explicit evaluation equals the Verlinde numbers. \vskip.1in

\subsection{The Verlinde numbers} To start, let $\quot$ be the Grothendieck Quot scheme parametrizing short exact sequences $$0
\rightarrow E \rightarrow \mathcal O^{r+k} \rightarrow F \rightarrow 0 \, \, \text{on} \, \, C,$$ with $E$ a subsheaf of rank $k$
and degree $-d.$ $\quot$ provides a compactification of the scheme $\mor$ of degree $d$ morphisms from $C$ to the Grassmannian $G(k, r+k)$. The morphism space sits inside $\quot$ as the locus of locally-free quotients of $\mathcal O^{r+k}.$ Throughout this
paper we will assume that the degree $d$ is sufficiently large to ensure that $\quot$ is irreducible \cite{BDW} and of the expected
dimension
$$\chi(E^{\vee}\otimes F)=(r+k)d-rk\gbar.$$ We will also assume that $d$ is divisible by $k$, and we set
$$s=(r+k)\cdot\frac{d}{k}-r\gbar.$$

As a fine moduli space, $\quot$ comes equipped with a universal sequence \begin{equation}\label{universal}0 \rightarrow \mathcal E
\rightarrow \mathcal O^{r+k} \rightarrow \mathcal F \rightarrow 0 \, \, \text{on} \, \, \quot \times C.\end{equation} We define the
cohomology class $$a_k = c_k ({\mathcal E^{\vee}|}_{\quot \times \{\text{point} \}}).$$ The intersection theory of the Chern classes of
$ {\mathcal E^{\vee}|}_{\quot \times \{\text{point} \}}$ has been well studied for over a decade: top intersections are computed by
the Vafa-Intriligator formula \cite{bertram}, \cite{sieberttian}, \cite{mo}, and have enumerative meaning
\cite{bertram},
counting degree $d$ maps from $C$ to the
Grassmannian $G(k, r+k)$ under incidence conditions with special Schubert subvarieties at fixed domain points.

Crucially for the arguments of this paper, the Verlinde numbers can be easily linked to the intersection theory of the Quot scheme.
This was
observed and justified through physical arguments by Witten \cite{witten}, who argued that the Wess-Zumino-Witten model for $GL(k)$
at level $r$ was related to the sigma model of the Grassmannian $G(k, r+k).$ More
precisely, we have the following spectacular realization of the Verlinde formula as a top intersection of the $a_k$ classes on   
$\quot$:

\begin{proposition}
\label{verl} For $d$ large enough, 
$$h^0 (SU (r, 0), \theta^k_r ) = \frac{r^g}{(k+r)^g} \int_{\quot} a_k^{s}.$$ \end {proposition}

{\it Proof.} We use the elementary expression for the Verlinde formula (the left-hand side of the above equation) appearing in
\cite {bl}, namely \begin{equation}\label{beauvilleverlinde}h^0 (SU (r, 0), \theta^k_r ) =
\frac{r^{g}}{(r+k)^{g}}\sum_{\stackrel{S\cup T=\{0, 1, \ldots, r+k-1\}}{|S|=k}} \prod_{\stackrel{s\in S}{t\in T}}
\left|2\sin\pi\frac{s-t}{r+k}\right|^{\gbar}.\end{equation} Now the Vafa-Intriligator formula \cite {mo} expresses
the intersection number on the Quot scheme as a sum over (ordered tuples of) roots $\lambda_1, \ldots, \lambda_k$ of unity 
of order $r+k$:
$$\int_{\quot}a_k^s=(-1)^{\gbar\binom{k}{2}}(r+k)^{k\gbar}\sum_{\lambda_1, \ldots, \lambda_k}
\prod_{i<j}\left(\left(\frac{\lambda_i}{\lambda_j}\right)^{\frac{1}{2}}-\left(\frac{\lambda_j}{\lambda_i}\right)^{\frac{1}{2}}\right)^{-2\gbar}=$$
$$=(r+k)^{k\gbar} \sum_{\stackrel{S\subset \{0, 1, \ldots,r+k-1\}}{|S|=k}} \prod_{\stackrel{s, t\in S}{s<t}}
\left(2\sin\pi\frac{s-t}{r+k}\right)^{-2\gbar}.$$ The proposition follows easily by repeatedly applying the following trigonometric
identity: \begin{equation} \label{trigid} \prod_{p=1}^{r+k-1}\left(2\sin\frac{p\pi}{r+k}\right) =r+k.\end{equation} \qed

Proposition \ref{verl} allows for the interpretation of the intersection number 
\begin{equation}\label{intonquot}q \stackrel{\rm def} {=} \int_{\quot} a_k^{s}\end{equation} itself, without prefactors, as a modified Verlinde number. Indeed, we have the following

\begin {corollary} Let $\det:U(r,0)\to \jac^0(C)$ denote the determinant morphism. For any two line bundles $A$ and $B$ of 
degree $\gbar$, 
\begin{equation}
\label{corol1}
h^{0}(U(r,0), \Theta_{r,A}^{k}\otimes {\det}^{\star} \Theta_{1, B})= \int_{\quot} a_k^{s}.
\end{equation}
\end {corollary}

{\it Proof.} Let $\tau: SU(r, \mathcal O)\times \jac^0(C)\to U(r,0)$ be the tensor product map $$(E, L)\to E\otimes L.$$ Note that by the 
see-saw principle, 
$$\tau^{\star}\Theta_{r, A}\cong\theta_r\boxtimes\Theta_{1, A}^{r}.$$ Indeed, it is easy to see that
the restriction of $\tau^{\star}\Theta_{r, A}$ to $SU(r, \mathcal O)\otimes \{L\}$ is isomorphic to $\theta_r$, while the restriction
$\tau^{\star}\Theta_{r, A}$ to $\{E\}\times \jac^{0}(C)$ is $$\Theta_{1, E\otimes A}\cong\Theta_{1,A}^r.$$ Therefore, letting
$$r:\jac^{0}(C)\to\jac^{0}(C)$$ denote multiplication by $r$, we have \begin{equation}\label{pullb}\tau^{\star}\left(\Theta_{r,A}^{k}\otimes
{\det}^{\star}\Theta_{1,B}\right)\cong\theta_r^{k}\boxtimes \Theta_{1,A}^{rk}\otimes r^{\star}\Theta_{1,B}.\end{equation} Since $\tau$ is an etale cover of
degree $r^{2g}$ we conclude $$\chi(\Theta_{r,A}^{k}\otimes {\det}^{\star} \Theta_{1, B})=\frac{1}{r^{2g}}\, \chi\left(\tau^{\star}
(\Theta_{r,A}^{k}\otimes {\det}^{\star}\Theta_{1,B})\right)=\frac{1}{r^{2g}} \,\chi(\Theta_{1,A}^{rk}\otimes
r^{\star}\Theta_{1,B})\, \chi(\theta_r^{k})=$$ $$= \frac{(rk+r^2)^g}{r^{2g}}\, \chi(\theta_r^{k})=\int_{\quot} a_k^{s}.$$

It remains to observe that the higher cohomology vanishes. By the Leray spectral sequence, this can be checked after pulling back by the finite
morphism $\tau$. The statement follows using the K{\"{u}}nneth decomposition, as the Theta bundles on $SU(r, \mathcal O)$ and the Jacobian have 
no
higher cohomologies. \qed

\begin{corollary}\label{sd2} For any line bundles $A, B, C, D$ of degree $\gbar$, we have $$h^{0}(U(r,0), \Theta_{r,A}^{k}\otimes {\det}^{\star} \Theta_{1, B})=h^{0}(U(k,0), \Theta_{k,C}^{r}\otimes {\det}^{\star} \Theta_{1, D}).$$
\end{corollary}

{\it Proof.}  This follows from \eqref{beauvilleverlinde} by interchanging $S$ and $T$.\qed\vskip.1in

\subsection {Transversality} We now set out to realize the intersection number \eqref{intonquot} geometrically. 
We will represent the cohomology class $a_k$ in different ways as vanishing loci of sections of the dual universal bundle 
$\mathcal E^{\vee}|_{\quot\times \{\text{point}\}}$. To this end, \begin{itemize} \item pick distinct points $p_1, \ldots, 
p_s$ on $C$. \item for each $1\leq i\leq s$, let $L_i$ be a one-dimensional subspace of $V=(\mathbb C^{r+k})^{\vee}$. This 
gives (by duality) a section $$\varphi_i:L_i\otimes \mathcal O\to V\otimes\mathcal O\to \mathcal E^{\vee}$$ on $\quot\times 
C$. \item consider the restricted morphisms $$\varphi_{i,p_i}:L_i\otimes \mathcal O\to \mathcal E^{\vee}_{p_i}$$ on $\quot 
\times \{p_i\}$, and further let $\mathcal Z_i$ denote the zero locus of the section $\varphi_{i,p_i}$ on $\quot\times 
\{p_i\}.$ \end {itemize}

As shown in \cite {bertram}, the schemes $\mathcal Z_i$ are of the expected codimension $k$, and
moreover $$\left[\mathcal Z_i\right]=a_k\cap \left[\quot\right].$$ We aim first to prove the following transversality statement, which is 
essentially contained in \cite{bertram}.

\begin{proposition} \label{inters0}For generic choices of subspaces $L_i$, $1\leq i\leq s$, the subschemes $\mathcal Z_1, \ldots, \mathcal Z_s$ 
representing 
the class $a_k$ intersect properly in reduced points, inside the locus of the Quot scheme consisting of stable vector bundle quotients.
\end{proposition}

{\it Proof.}  According to \cite{bertram}, 
generic zero-locus representatives for the classes $a_k$ appearing in the intersection on the right-hand side of 
\eqref{intonquot}, do not intersect outside the open subscheme $\mor$ of locally-free quotients of $\quot$.
We consider the restriction of each section $$\varphi_{i, p_i}: L_i \otimes \mathcal O\rightarrow \mathcal E^{\vee}_{p_i}$$ to $\mor$. For any point $p \in C,$ the evaluation of morphisms in the scheme $\mor$ gives a map
$$ev_p: \mor \rightarrow G(k,r+k).$$
The restriction of the zero locus $\mathcal Z_i$ of $\varphi_{i,p_i}$ to $\mor$ is 
$$\mathcal Z_i\big{|}_{\mor} = ev_{p_i}^{-1} (\mathcal W_i), $$
where $\mathcal W_i$ is the 
zero locus of the corresponding section $$\tilde{\varphi_i}: L_i \otimes \mathcal O\rightarrow V \otimes \mathcal O\rightarrow T^{\vee}$$ for the dual 
tautological bundle $T^{\vee}$
on $G(k, r+k).$

For each $i$, the zero locus $\mathcal W_i$ of $\tilde{\varphi_i}$ is a subGrassmannian $G(k, r+k-1)$ of $G(k, r+k)$. By Kleiman transversality
\cite{kleiman}, given a subscheme $Y$ of $\mor$, a generic subGrassmannian $G(k, r+k-1)$ intersects each locus in the image $ev_p (Y)$ over
which $ev_p$ has constant fiber dimension properly ({\it i.e.} in the correct codimension), hence $Y$ intersects $ev_p^{-1} (G(k, r+k-1))$
properly.  A top dimensional intersection of generic $\mathcal Z_i$s therefore avoids the unstable locus of $\mor$, and lies in the stable locus
$\mors$ consisting of short exact sequences \begin{equation} \label{ses} 0\to E\to \mathcal O^{r+k}\to F\to 0 \end{equation} with $E$ and $F$ 
stable.

Note now that for sufficiently large $d$, the restricted evaluation $$ev_p: \mors \rightarrow G(k, r+k)$$ has smooth fibers for every $p \in C.$
Indeed, the domain is smooth since the obstruction space is $$H^{1}(E^{\vee}\otimes F)=0.$$ This follows by considering the long exact sequence 
in cohomology associated with \eqref{ses} and noting that $H^{1}(F)=0$ when $F$ is stable and $d$ is large enough. The conclusion about $ev_p$ 
is implied by generic smoothness and
by the fact that the target is a homogeneous space.

Now, if $Y$ is a generically smooth subscheme of $\mors$, then $ev_p (Y)$ is generically smooth in $G(k, r+k)$, and 
by Kleiman transversality we conclude that the intersection of $ev_p(Y)$ with a generic subGrassmannian $G(k, r+k-1)$ is generically smooth,
implying by the smoothness of the fibers of $ev_p$ that the intersection of $Y$ with $ev_p^{-1} (G (k, r+k-1))$ is generically smooth as well. 
Inductively then, the intersection of generic 
$\mathcal Z_i, \, 1\leq i \leq s$ is a zero-dimensional generically smooth subscheme of $\mors$. Therefore, this intersection is smooth. 
This establishes Proposition \ref{inters0}. \qed
\vskip.1in

Considering all the sections $$\varphi_{i, p_i}: L_i\otimes \mathcal O\rightarrow \mathcal E^{\vee}_{p_i}, \, \, 1\leq i\leq s$$ at once, we 
obtain that $$\mathcal Z=\cap_{i=1}^{s}\mathcal Z_i$$ is the smooth zero locus on 
$\mors$ of a section
$$\varphi: \mathcal O\rightarrow \mathcal E^{\vee}_{p_1} \oplus \cdots \oplus \mathcal E^{\vee}_{p_s}.$$ Therefore, we have an exact sequence on $\mathcal Z$
\begin{equation}
\label{tangent1}
0\rightarrow T\mathcal Z \rightarrow T\mors\big{|}_{\mathcal Z} \rightarrow \mathcal E^{\vee}_{p_1} \oplus \cdots \oplus \mathcal E^{\vee}_{p_s}\big{|}_{\mathcal Z} \rightarrow 0.
\end{equation}
We would like to identify $T\mathcal Z$ in terms of the available universal 
structures. 

Each of the sections $\varphi_i: L_i \otimes \mathcal O\rightarrow \mathcal E^{\vee}$ on $\quot \times C$ vanishes on 
$\mathcal Z\times \{p_i\}$. On $\mathcal Z \times C$ therefore, each $\varphi_i$ factors as $$L_i \otimes \mathcal 
O\rightarrow L_i \otimes {\pi}_{C}^{\star}\mathcal O(p_i) \rightarrow \mathcal E^{\vee},$$ where ${\pi}_{C}$ denotes the 
projection $$\pi_C:\mathcal Z\times C \rightarrow C.$$ Let $S$ denote the pushout on $C$ of the two maps 
$$\bigoplus_{i=1}^{s}L_i\otimes \mathcal O\to \bigoplus_{i=1}^{s}L_i\otimes\mathcal O(p_i)$$ and 
$$\bigoplus_{i=1}^{s}L_i\otimes \mathcal O\to V\otimes \mathcal O.$$ Note that \begin{equation} \det S = \mathcal O(p_1 + 
\cdots + p_s), \end{equation} and that for distinct $p_i$s, $S$ is locally free. Considering the $s$ 
conditions imposed by each one of the points $p_i$, we arrive at the 
conclusion that the universal map $$V\otimes \mathcal O\to \mathcal E^{\vee}\to 0$$ factors as $$V\otimes \mathcal O\to 
\pi_{C}^\star S\to \mathcal E^{\vee}\to 0 \text{ on }\mathcal Z\times C.$$ Let ${\mathcal F'}$ on $\mathcal Z \times C$ 
denote the ensuing kernel,\begin{equation}\label{kernel} 0 \rightarrow {\mathcal F'} \rightarrow {\pi}_{C}^{\star} S 
\rightarrow \mathcal E^{\vee}\to 0.\end{equation}

Denoting by $\iota: S^{\vee} \hookrightarrow \mathcal O^{r+k}$ the 
induced inclusion of sheaves, we see that each exact sequence $$0 \rightarrow E \rightarrow \mathcal O^{r+k} \rightarrow F \rightarrow 0$$ representing 
a point $\zeta$ in the 
intersection $\mathcal Z = \bigcap_{i=1}^{s}\mathcal Z_i$ factors through $\iota$. For each $\zeta \in Z$ there is thus an associated exact 
sequence $$0 \rightarrow E \rightarrow S^{\vee} \rightarrow F'^{\vee} \rightarrow 0.$$ 

\begin{lemma}
\label{zariski}

The Zariski tangent space to $\mathcal Z$ at $\zeta$ is 
$$T_{\zeta}\mathcal Z= H^0 (E^{\vee} \otimes F'^{\vee}).$$
\end{lemma}

{\it Proof.}  As a consequence of \eqref{universal}, \eqref{kernel} and the five lemma, we obtain a natural exact sequence
\begin{equation}
\label{newf}
0 \rightarrow \mathcal F'^{\vee} \rightarrow \mathcal F\rightarrow \pi_{C}^{\star} \mathcal O_{p_1 + \ldots + p_s} \rightarrow 0 \, \, \, \text{on} \, \, \mathcal Z 
\times C,
\end{equation}
where the second morphism is given by pairing evaluations of sections of $\mathcal F$ at the points $p_1, \ldots p_s$ with $\varphi_1, \ldots , \varphi_s$ respectively.
Twisting the exact sequence \eqref{newf} with $\mathcal E^{\vee}$ and pushing forward to $\mathcal Z$ by the projection $\pi_{\mathcal Z}: \mathcal Z \times C \rightarrow \mathcal Z$,
we get
\begin{equation}
\label{push}
0 \rightarrow R^0{\pi_{\mathcal Z}}_{\star} (\mathcal E^{\vee} \otimes \mathcal F'^{\vee}) \rightarrow  R^0{\pi_{\mathcal Z}}_{\star} (\mathcal E^{\vee} \otimes \mathcal F) \rightarrow 
\oplus_{i=1}^{s} \mathcal E^{\vee}_{p_i} \rightarrow \cdots
\end{equation}
Comparing \eqref{tangent1} and \eqref{push} we conclude that the sequence \eqref{push} is exact on the right and that 
$$T\mathcal Z = R^0{\pi_{\mathcal Z}}_{\star} (\mathcal E^{\vee} \otimes \mathcal F'^{\vee}).$$\qed\vskip.1in

The smoothness and zero-dimensionality of ${\mathcal Z} = {\mathcal Z}_1 \cap \ldots \cap {\mathcal Z}_s$ established by Proposition 
\ref{inters0} implies that in Lemma \ref{zariski} we have
\begin{equation}
\label{zerots}
T_{\zeta}\mathcal Z = H^0 (E^{\vee} \otimes F'^{\vee}) = 0.
\end{equation}
After dropping the primes and duals on the $F$s, our discussion ensures the existence of $q$ 
exact sequences 
\begin{equation}
\label{rightbun}
0 \rightarrow E_i\rightarrow S^{\vee}
\rightarrow F_i \rightarrow 0, \, \, \, 1 \leq i \leq q,
\end{equation}
 with \begin{itemize} \item[(i)] $E_i$ is a {\it stable} bundle of degree $-d$ and rank
$k$, \item[(ii)] $F_i$ is a bundle of degree $d-s$ and rank $r$, \item[(iii)] $\det E_i\otimes \det F_i=\det S^{\vee}=\mathcal O(-p_1-\ldots-p_s),$ \item[(iv)]
$ h^0(E_i^{\vee} \otimes F_j)= 0 \, \, \text{iff} \, \, i=j. $ \end{itemize}
The last assertion follows from \eqref{zerots} and the observation that for $i\neq j$ we have a non-zero map $$E_i\to S^{\vee}\to F_j.$$
Property $\text{(iv)}$ also proves the stability of the bundle $F_i$ via lemma $A.1$ in \cite {Bel}, as the required 
hypothesis $h^{0}(E_i^{\vee}\otimes
F_i)=h^1(E_i^{\vee}\otimes F_i)=0$ is satisfied. \vskip.1in

\begin{remark} The tangent space computation of Lemma \ref{zariski} and the arguments preceding it identify $\mathcal Z$ as 
the Quot scheme ${\text{Quot}}_{d} (S^{\vee}, k, C)$ of rank $k$ and degree $-d$ subsheaves of $S^{\vee}$. The pairs $(E_i, 
F_i)$ constructed above correspond to the closed points of this smooth zero dimensional Quot scheme of length equal to the 
Verlinde number \eqref{corol1}. \end{remark}

Finally, we construct the pairs of bundles that were used in the proof of Theorem \ref{main2}.

\begin{proposition}\label{pairs}
There exist $q$ pairs of bundles $(A_1, B_1), \ldots, (A_q, B_q)$ with the following properties
\begin{itemize} \item $A_1, \ldots, A_q$ are rank $r$ degree $0$ stable bundles 
\item $B_1, \ldots, B_q$ are rank $k$ degree $0$ stable bundles 
\item $h^{0}(A_i\otimes B_j\otimes M)=0$ iff $i=j$
\item $h^{0}(\det A_i^{\vee}\otimes \det B_i\otimes L^{-1}\otimes K)=0.$
\end {itemize}
\end {proposition}

{\it Proof. } We only need to adjust the degrees of the bundles \eqref{rightbun}. Pick a line bundle $Q$ of degree 
$-\frac{d}{k}$. This ensures that the line bundle $$R=M^{-r}\otimes Q^{-r-k}\otimes L^{-1}\otimes K(-p_1-\ldots-p_s)$$ has 
degree $\gbar$. Moreover, we may assume that $Q$ is chosen so that $R$ avoids the theta divisor on $\jac^{\gbar}(C)$ 
{\it{i.e.}}, so that it has no global sections. Setting $$A_i=E_i^{\vee}\otimes Q, \text { and } B_i= F_i\otimes 
Q^{-1}\otimes 
M^{-1}$$ we obtain, via Property (iv) above, that $$h^{0}(A_i\otimes B_j\otimes M)=h^{0}(E_i^{\vee}\otimes F_j)=0\text{ if 
and only if } i=j.$$ Moreover $$\det A_i^{\vee}\otimes \det B_i=\det E_i \otimes \det F_i\otimes M^{-r}\otimes Q^{-r-k} 
=\det S^{\vee}\otimes M^{-r}\otimes Q^{-r-k}.$$ Therefore, $$h^{0}(\det A_i^{\vee}\otimes \det B_i\otimes L^{-1}\otimes 
K)=h^{0}( M^{-r}\otimes Q^{-r-k}\otimes L^{-1}\otimes K(-p_1-\ldots-p_s))=h^{0}(R)=0.$$ This is what we set out to prove. 
\qed

\section{The argument on the Jacobian} In this section we prove Theorem \ref{main}. We pass to the moduli spaces of 
semistable bundles with fixed determinant via the natural inclusion $$j:SU(r,\mathcal O)\hookrightarrow U(r, 0).$$ 
Restricting sections gives a morphism (well defined projectively) $$\rho : H^{0}(U(r, 0), \Theta_{r, M}^k \otimes 
{\det}^{\star} \Theta_{1, L} ) \rightarrow H^0 (SU (r, \mathcal O), \theta^k_r ).$$ We show

\begin {proposition}\label{inj} For generic $L$ and $M$ the restriction $\rho$ is surjective. 
\end {proposition}

{\it Proof.} Set $$\mathcal L= \Theta_{r, M}^k \otimes
{\det}^{\star} \Theta_{1, L}.$$ Recall from \eqref{pullb} that the pullback of $\mathcal L$ under the tensor product map  
$$\tau: SU (r, \mathcal O) \times \jac^{0}(C) \rightarrow U (r, 0)$$ is given by 
\begin{equation}
\label{isomchoice}
\tau^{\star} \mathcal L\cong \theta^k_r \boxtimes \Theta^{kr}_{1,M} \otimes 
r^{\star}\Theta_{1,L}.
\end{equation}
Such an isomorphism will be fixed below. 

Let $\mathsf G$ be the group of $r$-torsion points in $\jac^0(C)$. Under the natural antidiagonal action of $\mathsf G$ on 
$SU (r, \mathcal O) 
\times \jac^{0}(C)$, the morphism $\tau$ is a Galois covering with Galois group $\mathsf G$.
It follows that $$H^{0}(U(r,0), \mathcal L)=H^{0}(SU(r, \mathcal O)\times \jac^{0}(C), \tau^{\star}\mathcal L)^{\mathsf G},$$ the right hand
side being the space of $\mathsf G$-invariant sections. It suffices to show that given any section $s$ of $\theta_r^k$ on $SU (r, \mathcal O)$, there is a
$\mathsf G$-invariant section $\tilde{s}$ of $\tau^{\star} \mathcal L$ which restricts to $s$ over $SU (r, \mathcal O) \times \{\mathcal O\}.$

First, we claim that for a generic choice of $L$, we can find a section $u$ of $r^{\star}\Theta_{1,L}$ which does not 
vanish at the origin, but vanishes at all other points of $\mathsf G$. This claim follows if we show that the sheaf 
$r^{\star}\Theta_{1, L} \otimes {\mathcal O}_{\mathsf G}$ is generated by the global sections of $r^{\star}\Theta_{1, L}.$ 
Equivalently we prove that the map $$\mu: H^0 (\jac^0(C), r^{\star}\Theta_{1,L}) \rightarrow H^{0}(\jac^{0}(C), 
r^{\star}\Theta_{1, L}\otimes \mathcal O_{\mathsf G})$$ is surjective. Since both $H^0 (\jac^{0}(C), 
r^{\star}\Theta_{1,L})$ and $H^{0}(\jac^{0}(C), r^{\star}\Theta_{1, L}\otimes \mathcal O_{\mathsf G})$ have dimension 
$r^{2g}$, it suffices to explain that the morphism $\mu$ is injective. This is a consequence of the vanishing 
\begin{equation}\label{vanish}H^{0}(\jac^{0}(C), r^{\star}\Theta_{1,L}\otimes I_{\mathsf G})=0,\end{equation} combined with 
the long exact sequence in cohomology associated with \begin{equation} \label{ideals} 0\to I_{\mathsf G}\to \mathcal 
O_{\jac} \to \mathcal O_{\mathsf G}\to 0. \end{equation} Here, $I_{\mathsf G}$ is the ideal sheaf of $\mathsf
G$ in ${\text{Jac}}^0 (C).$
To prove \eqref{vanish}, we observe that since $r$ is a finite morphism $$H^{0}(\jac^0(C), 
r^{\star}\Theta_{1,L}\otimes I_{\mathsf G})=H^{0}(\jac^0(C), r_{\star}\left(r^{\star}\Theta_{1,L}\otimes I_{\mathsf 
G}\right)).$$ Now it is well known that $$r_{\star} \mathcal O_{\jac}=\bigoplus_{i} T_i$$ where $T_{i}$ are the $r$-torsion 
line bundles on $\jac^{0}(C)$ \cite {DT}, under the standard identification of degree $0$ line bundles on the Jacobian and 
on the curve.  From the pushforward under $r$ of the ideal sheaf sequence \eqref{ideals}, it easily follows that 
$$r_{\star} I_{\mathsf G}=\bigoplus_{i} T_i\otimes I_{\{\mathcal O\}},$$ with $I_{\{\mathcal O\}}$ being the ideal sheaf of 
the origin $\{\mathcal O\}\in \jac^0 (C)$. Hence, $$H^{0}(\jac, r^{\star}\Theta_{1,L}\otimes I_{\mathsf G})=\bigoplus_{i} 
H^{0}(\jac,\Theta_{1, L}\otimes T_i\otimes I_{\{\mathcal O\}})=\bigoplus_{i}H^{0}(\jac, \Theta_{1, L\otimes T_i}\otimes 
I_{\{\mathcal{O}\}})=0.$$ The last equality follows provided the (canonical) section of $\Theta_{1,L\otimes T_i}$ does not 
vanish at the origin $\{\mathcal O\}$. The nonvanishing is equivalent to $h^{0}(L\otimes T_i)=0$, which can be 
arranged if we pick $L$ outside the Theta divisors $\Theta_{1, T_i}$ for all $i$. This proves the existence of the section 
$u$.

Next, take a section $v$ of $\Theta_{1,M}^{kr}$ which does not 
vanish at the origin. Again, this can be arranged for generic $M$. The section $t=u\otimes v$ of 
$r^{\star}\Theta_{1,L}\otimes \Theta_{1,M}^{kr}$ does not vanish at the origin, but vanishes at all the other points in $\mathsf G$. 

Finally, pick any section $s$ of the line bundle $\theta_{r}^{k}$ on $SU(r, \mathcal O)$. We extend $s$ to a $\mathsf 
G$-invariant section $\tilde s$ by averaging over $\mathsf G$. First, consider the section $s\boxtimes t$ of the line 
bundle $$\theta_{r}^k\boxtimes r^{\star}\Theta_{1,L}\otimes \Theta_{1, M}^{kr}\cong\tau^{\star} \mathcal L.$$ For each 
$g\in \mathsf G$, the pullback $g^{\star}(s\boxtimes t)$ is a section of $g^{\star} \tau^{\star} \mathcal L=\tau^{\star} 
\mathcal L$. The average $$\tilde{s}=\sum_{g\in \mathsf G} g^{\star}(s\boxtimes t)$$ is then an invariant section of 
$\tau^{\star} \mathcal L$. By construction, the pullback section $g^{\star} (s \boxtimes t)$ vanishes along $SU(r, \mathcal 
O)\times \{\mathcal O\}$ if and only if $g \in {\mathsf G} \setminus \{\text{identity}\}$. Thus the two sections 
$j^{\star}\tilde{s}$ and $s$ differ by a nonzero constant, as wished. \qed\vskip.1in

As a consequence of the lemma, the dual restriction map
$$\sigma: {H^0 (SU (r, {\mathcal O}), \theta_r^k )}^{\vee} \rightarrow {H^0 (U(r, 0), \Theta_{r, M}^k \otimes 
{\det}^{\star}\Theta_{1, L})}^{\vee}$$ is injective. This is important for the proof of Theorem \ref{main}, which we now give. \vskip.1in

{\it Proof of Theorem \ref{main}.} Consider the injection $$\iota:H^0 (U (k, 0), \Theta_{k,M}^r)\hookrightarrow H^0 (U (k, 
0), \Theta_{k,M}^r \otimes {\det}^{\star} \Theta_{1, L^{-1}\otimes K} )$$ given by $$s \mapsto s \otimes {\det}^{\star} 
\Theta_{1, L^{-1}\otimes K},$$ where $\Theta_{1, L^{-1}\otimes K}$ now denotes the canonical section of the Theta line 
bundle. Note that via the identification $U(k, 0)\cong U(k, k\gbar)$ given by tensoring bundles with $M$, we have 
\begin{equation} \label{equality}\mathsf {D}\circ \sigma=\iota\circ\mathsf {SD}.\end{equation} Here $\mathsf {SD}$ is the 
strange duality morphism \eqref{sd}. 
The above equality follows 
directly from the definitions. Indeed, we observe that the restriction of the divisor $\Delta$ to $SU(r, \mathcal O)\times 
U(k,0)$, which determines the morphism $\mathsf{D}\circ \sigma$, equals $$\{(E,F) \text { such that } h^{0}(E\otimes 
F\otimes M)\neq 0 \text { or } h^{0}(\det F\otimes L^{-1}\otimes K)\neq 0\}.$$ But this is also the divisor inducing the 
morphism $\iota\circ \mathsf{SD}$.

Finally, by Theorem \ref{main2} and Proposition \ref{inj}, the composition $$\mathsf {D} \circ \sigma: H^0 (SU (r, \mathcal O), \theta^k_r )^{\vee} \rightarrow H^0 (U (k, 0),
\Theta_{k,M}^r\otimes {\det}^{\star}\Theta_{1,L^{-1}\otimes K})$$ is an isomorphism onto its image. Therefore, \eqref{equality} implies that $\mathsf {SD}$ is injective. A well known symmetry of the Verlinde
formula gives the equality of dimensions \cite {Bea} $$h^{0}(SU(r,\mathcal O), \theta_{r}^{k})=h^{0}(U(k,0), \Theta_{k, M}^{r})\cong h^{0}(U(k, k\gbar), \Theta_k^r).$$ Hence $\mathsf{SD}$ must be an isomorphism. This proves Theorem \ref{main}.\qed

\section{Strange duality in arbitrary degree} 

The arguments of this paper can be used to establish the strange duality conjecture for arbitrary rank and degree. In this 
section,
we will indicate the relevant statements and their proofs.

Consider $d$ and $r$ two relatively prime integers, and $h, k$ two non-negative integers. Without loss of generality assume  
that 
$0<d<r$. To start, we define an extended duality morphism involving the moduli spaces of bundles $U(hr, hd)$ and $U(kr, 
-kd)$. Fix a 
reference bundle $Q$ of degree $d$ and rank $r$, and fix two line bundles $L$ and $M$ of degree $\gbar$. The duality 
relates sections of the pluri-Theta bundles $$\Theta_{hr, M\otimes Q^{\vee}}^{k}\to U(hr, hd), \text{ and }\Theta_{kr, 
M\otimes Q}^{h}\to U(kr, -kd),$$ twisted by line bundles coming from the Jacobian via the determinant maps $$\det: U(hr, 
hd)\to \jac^{hd}(C)\text { and }\det: U(kr, -kd)\to \jac^{-kd}(C).$$

Consider the divisor \begin{eqnarray}\Delta=\{(E, F)&\in& U(hr, hd)\times U(kr, -kd) \text { such that } h^{0}(E\otimes 
F\otimes M)\neq 0 \text { or }\nonumber\\&&h^{0}(\det E^{\vee}\otimes \det F\otimes (\det Q)^{h+k}\otimes L^{-1}\otimes 
K)\neq 0\}.\nonumber\end{eqnarray} It is easy to check, in complete analogy with the argument of Lemma \ref{cross}, that 
$$\mathcal O(\Delta)=\left(\Theta_{hr, M\otimes 
Q^{\vee}}^{k}\otimes {\det}^{\star} \Theta_{1, L\otimes {\left(\det Q\right)}^{-h}}\right)\boxtimes\left(\Theta_{kr, 
M\otimes Q}^{h}\otimes {\det}^{\star}\Theta_{1, L^{-1}\otimes K \otimes \left(\det Q\right)^{k}}\right).$$ Note that we 
recover the 
divisor used in the previous sections, when the degree is $0$, by taking $Q$ to be trivial. \vskip.1in

We will show

\addtocounter{theorem}{-2}\begin{theorem}\label{main3} The divisor $\Delta$ induces a rank-level duality isomorphism between
$$H^{0}\left(U(hr, hd), \Theta_{hr, M\otimes Q^{\vee}}^{k}\otimes {\det}^{\star} \Theta_{1, L\otimes {\left(\det
Q\right)}^{-h}}\right)^{\vee}, \text { and }$$ $$H^{0}\left(U(kr, -kd), \Theta_{kr, M\otimes Q}^{h}\otimes
{\det}^{\star}\Theta_{1, L^{-1}\otimes K \otimes \left(\det Q\right)^{k}}\right).$$ \end{theorem}

The starting point of the argument is the following rewriting of the Verlinde formula, valid in arbitrary degree.
\begin{proposition}\label{verl1} Let $\theta_{hr}$ denote the ample generator for the Picard group of $SU(hr, hd)$. If $e\equiv -d 
\mod r$ is sufficiently large, then $$h^{0}(SU(hr, hd), \theta^k_{hr})=\frac{h^g}{(h+k)^g}\int_{\text{Quot}_{ke}(\mathcal
O^{r(h+k)}, kr, C)} a_{kr}^{(h+k)e-hr\gbar}.$$ \end {proposition}

{\it Proof.} Theorem $9.4$ in \cite {bl} establishes an isomorphism between the space of global sections and a space of conformal blocks for the affine Lie algebra $\widehat{\mathfrak{sl}}_{hr}$ at level $kr$. The dimension of the space of conformal blocks is computed in corollary $9.8$ of \cite {Bea2}. Together, these two statements give the formula
$$(r(h+k))^{hr\gbar}\left(\frac{h}{h+k}\right)^{\gbar}\sum_{\vec\nu}\left\{\text{Trace}_{\lambda}\left(\exp 2\pi i
\frac{\vec\nu}{r(h+k)}\right)\cdot\prod_{i<j}\left(2\sin \pi\frac{\nu_i-\nu_j}{r(h+k)}\right)^{-2\gbar}\right\}.$$ In the above   
summation $\vec{\nu}$ is an element of the weight lattice of $\mathfrak{sl}_{hr}$ whose standard coordinates satisfy the property
$$\nu_1>\ldots> \nu_{hr},\, \nu_i-\nu_j \in \mathbb Z,\,\nu_1-\nu_{hr}<r(h+k), \text { and } \nu_1+\ldots+\nu_{hr}=0.$$ Also, the 
trace is computed in the representation of $\mathfrak{sl}_{hr}$ whose highest weight is the following combination of the standard
coordinate forms $$\lambda=kr\left(e_1+\ldots+e_{hr-hd}\right).$$ The Weyl character formula expresses the trace
of $H=\exp\left(2\pi i\frac{\vec{\nu}}{r(h+k)}\right)$ as $$\text{Trace}_{\lambda}(H)=\frac{\sum_{\sigma\in
S_{hr}}\text{sgn}(\sigma)\cdot(\lambda+\rho)(\sigma H)}{\sum_{\sigma\in S_{hr}}\text{sgn}({\sigma})\cdot \rho(\sigma H)},$$ with
$\rho$ being half the sum of the positive roots of $\mathfrak{sl}_{hr}$. To simplify this formula, we compare the terms indexed by $\sigma\tau$ in the numerator with
terms indexed by $\sigma$ in the denominator, where $\tau$ is the permutation $$\tau(i)=\begin{cases}i+hd,
\text{ if } i\leq hr-hd,\\ i-hr+hd, \text{ if } i>hr-hd.\end{cases}.$$ The explicit computation shows that all such terms are proportional by
the factor $$(-1)^{hd(r-1)}\exp\left(2\pi i (\nu_1+\ldots+\nu_{hr-hd})\right).$$ We now apply the following change of 
variables, which can also be found in \cite {zagier}. Setting $$t_{i}=\nu_{i}-\nu_{hr}, \,\,1\leq i\leq hr,$$ the Verlinde formula becomes
$$(-1)^{hd(r-1)}(r(h+k))^{\gbar}\left(\frac{h}{h+k}\right)^{\gbar}\sum_{\stackrel{T\subset\{0, \ldots, r(h+k)-1\}}{0\in T,\,|T|=hr}}\exp\left(2\pi i\frac{d}{r}\sum_{t\in T} t\right)\prod_{i<j}\left(2\sin\pi \frac{t_i-t_j}{r(h+k)}\right)^{-2\gbar}.$$ As in \cite {zagier}, we
remove the condition $0\in T$ at the expense of including the prefactor $\frac{hr}{r(h+k)}$. Further, we rewrite the expression in a form which is symmetric in $T$ and its complement $S$, using the trigonometric identity \eqref{trigid}. We arrive at the
equality \begin{eqnarray}\label{ver2} h^{0}(SU(hr, hd), \theta^k_{hr})&=&(-1)^{hd(r-1)}\frac{h^g}{(h+k)^{g}}\sum_{\stackrel{S\cup
T=\{0, \ldots, r(h+k)-1\}}{|T|=hr}} \exp\left(2\pi i \frac{d}{r}\sum_{t\in T} t\right)\times\nonumber\\&\times&\prod_{\stackrel{s\in
S}{t\in T}} \left|2\sin\pi\frac{s-t}{r(h+k)}\right|^{\gbar}.\end{eqnarray} Finally, the Lemma follows from \eqref{ver2} by a 
backward use of the Vafa-Intriligator formula, just as in Proposition \ref{verl}.\qed\vskip.1in

{\it Proof of Theorem \ref{main3}. } The proof is similar to that of Theorem \ref{main2}, but the numerical details are slightly
different. By Proposition \ref {verl1}, and the usual covering argument, we have $$h^{0}\left(U(hr, hd), \Theta_{hr, M\otimes 
Q^{\vee}}^{k}\otimes {\det}^{\star} \Theta_{1, L\otimes {\left(\det Q\right)}^{-h}}\right)=\int_{\text{Quot}_{ke}(\mathcal
O^{r(h+k)}, kr, C)} a_{kr}^{(h+k)e-hr\gbar}.$$ It follows from \eqref{ver2}, after interchanging $S$ and $T$, that
$$h^{0}\left(U(kr, -kd), \Theta_{kr, M\otimes Q}^{h}\otimes {\det}^{\star}\Theta_{1, L^{-1}\otimes K \otimes \left(\det
Q\right)^{k}}\right)$$ is computed by the {\it same} intersection number on Quot. Let $q$ denote this common value.

Finally, the reduced points representing the intersection number $q$ give, after adjusting the degrees, $q$ pairs of 
stable bundles $$(A_i, B_i)\in U(hr, hd)\times U(kr, -kd), 1\leq i\leq q,$$ such that $\Delta_{A_i}(B_j) \neq 0$ if and 
only if $i = j.$ This is
enough to complete the proof. \qed\vskip.1in

Theorem \ref{main3} leads to the proof of the arbitrary degree strange duality conjecture, as stated in \cite {DT}. As usual, let 
$$\tau:SU\left(hr, (\det Q)^{h}\right)\times U(kr, k(r\bar{g} - d))\to U(khr^2, khr^2\bar{g})$$ be the tensor product map
$$(E,F)\to E\otimes F.$$ Using the see-saw theorem, we compute $$\tau^{\star} \Theta_{khr^2}\cong\theta^{k}_{hr}\boxtimes
\Theta_{kr, Q}^h.$$ The canonical pullback section of $\tau^{\star}\Theta_{khr^2}$ induces a morphism between the spaces of
global sections: $$\mathsf{SD}:H^{0}\left(SU\left(hr, (\det Q)^{h}\right), \theta^k_{hr}\right)^{\vee}\to H^{0}\left(U\left(kr,
k(r\bar{g} - d)\right), \Theta_{kr,Q}^{h}\right).$$

\begin{theorem} $\mathsf {SD}$ is an isomorphism.
\end {theorem}

{\it Proof. } Restrict the duality morphism of Theorem \ref{main3} to $$j:SU(hr, (\det Q)^{h})\hookrightarrow U(hr, hd)$$ and
repeat the argument of Section $4$. \qed


\begin{thebibliography}{1}

\bibitem [Bea1]{Bea1}

A. Beauville, {\it Fibr{\'{e}}s de rang $2$ sur les courbes, fibr{\'{e}} d{\'{e}}terminant et fonctions 
th{\^{e}}ta}, Bull. Soc. Math. France 116 (1988), no. 4, 431-448. 

\bibitem [Bea2]{Bea}

A. Beauville, {\it Vector bundles on curves and generalized theta functions: recent results and open problems}, Current topics in complex algebraic geometry, 17-33, Math. Sci. Res. Inst. Publ., 28, Cambridge Univ. Press, Cambridge, 1995.

\bibitem [Bea3]{Bea2}

A. Beauville, {\it Conformal blocks, fusion rules and the Verlinde formula}, Proceedings of the Hirzebruch 65 Conference on Algebraic Geometry, 75-96, Israel Math. Conf. Proc., 9, Bar-Ilan Univ., Ramat Gan, 1996.

\bibitem [BL]{bl}

A. Beauville, Y. Laszlo, {\it Conformal blocks and generalized theta functions}, Comm. Math. Phys. 164 (1994), no. 2, 385-419.

\bibitem [BNR]{BNR}

A. Beauville, M.S. Narasimhan, S. Ramanan, {\it Spectral curves and the generalised theta divisor}, J. Reine Angew. Math. 398 (1989), 169-179.

\bibitem [Bel]{Bel}

P. Belkale, {\it The strange duality conjecture for generic curves}, AG/0602018.

\bibitem [Ber]{bertram} A. Bertram, {\it Towards a Schubert Calculus for Maps
from a Riemann Surface to a Grassmannian}, {Internat. J. Math {5}} {(1994)},
{no 6}, {811-825}.

\bibitem [BDW]{BDW}
A. Bertram, G. Daskalopoulos and R. Wentworth, {\it Gromov Invariants for
Holomorphic Maps from Riemann Surfaces to Grassmannians}, {J. Amer. Math. Soc. {{9}}}
{(1996)}, {no 2}, {529-571}.

\bibitem [DN]{DN}

J. M. Drezet, M. S. Narasimhan, {\it Groupe de Picard des vari{\'{e}}t{\'{e}}s de modules de fibr{\'{e}}s 
semi-stables sur les courbes alg{\'{e}}briques}, Invent. Math. 97 (1989), no. 1, 53-94. 

\bibitem [DT]{DT}

R. Donagi, L. Tu, {\it Theta functions for ${\rm SL}(n)$ versus ${\rm GL}(n)$}, Math. Res. Lett. 1 (1994), no. 3, 345-357.

\bibitem [F]{F}

I. Frenkel, {\it Representations of affine Lie algebras, Hecke modular forms and Korteweg-de Vries type equations}, Lecture Notes in Math., 933, Springer, Berlin-New York, 1982.

\bibitem [vGP]{vGP}

B. van Geemen, E. Previato, {\it Prym varieties and the Verlinde formula}, Math. Ann. 294 (1992), no. 4, 741-754.

\bibitem [K]{kleiman}

S. Kleiman, {\it The transversality of a general translate}, Compositio Math. 28 (1974), 287-297.

\bibitem [LB]{LB}

H. Lange, C. Birkenhake, {\it Complex abelian varieties}, Springer-Verlag, Berlin-New York, 1992. 

\bibitem [L]{l}

Y. Laszlo, {\it A propos de l'espace des modules de fibr{\'{e}}s de rang $2$ sur une courbe}, Math. Ann. 299 
(1994), no. 4, 597-608. 

\bibitem [MO]{mo}

A. Marian, D. Oprea, {\it Virtual intersections on the Quot scheme and Vafa-Intriligator formulas}, to be published in Duke Mathematical Journal, arxiv AG/0505685.

\bibitem [M]{M}

D. Mumford, {\it Tata lectures on theta I}, Progress in Mathematics, Birkhauser, Boston, 1983.
 
\bibitem [NS]{NS}

S. Naculich, H. Schnitzer, {\it Duality relations between ${\rm SU}(N)\sb k$ and ${\rm SU}(k)\sb N$ WZW models and their braid matrices}, Phys. Lett. B 244 (1990), no. 2, 235-240.

\bibitem [P]{P}

M. Popa, {\it Verlinde bundles and generalized Theta series}, Trans. Amer. Math. Soc, 254 (2002), 1869-1898.

\bibitem [RT]{RT}

B. Russo, M. Teixidor i Bigas, {\it On a conjecture of Lange}, J. Algebraic Geometry, 8 (1999), 483-496. 

\bibitem [ST]{sieberttian}
B. Siebert, G. Tian, {\it On quantum cohomology rings of Fano manifolds and
a formula of Vafa and Intriligator,} {Asian J. Math. {1}} {(1997)}, {no. 4},
{679--695.}

\bibitem [W]{witten}

E. Witten, {\it The Verlinde algebra and the cohomology of the Grassmannian}, Geometry, topology and physics, 357 - 422, 
Conf. Proc. Lecture Notes Geom. Topology, IV,  Internat. Press, Cambridge, MA, 1995.

\bibitem [Z] {zagier}

D. Zagier, {\it Elementary aspects of the Verlinde formula and of the Harder-Narasimhan-Atiyah-Bott formula,} 
Proceedings of the Hirzebruch 65 Conference on Algebraic Geometry, 445-462, Israel Math. Conf. Proc., 9, Bar-Ilan Univ., Ramat Gan, 1996.

\end{thebibliography}
\end{document}